\documentclass[journal]{IEEEtran}
\ifCLASSINFOpdf
\else
\fi

\usepackage{balance}
\usepackage{float}
\usepackage{caption}
\usepackage{subcaption}
\usepackage{tikz,pgfplots,pgfplotstable,booktabs}
\usepackage[official]{eurosym}
\usepackage{amsmath}
\usepackage{xcolor}
\usepackage{colortbl}
\usepackage{comment}
\allowdisplaybreaks[1]

\pgfplotsset{
    discard if not/.style 2 args={
        x filter/.code={
            \edef\tempa{\thisrow{#1}}
            \edef\tempb{#2}
            \ifx\tempa\tempb
            \else
                
            \fi
        }}}

\newenvironment{ldescription}[1]
  {\begin{list}{}%
   {\renewcommand\makelabel[1]{##1\hfill}%
   \settowidth\labelwidth{\makelabel{#1}}%
   \setlength\leftmargin{\labelwidth}
   \addtolength\leftmargin{\labelsep}}}
  {\end{list}}
  
\begin{document}
%
\title{Time-Adaptive Unit Commitment}
%
%
%

\author{Salvador Pineda, Ricardo Fern{\'a}ndez-Blanco, and
        Juan Miguel Morales
\thanks{S. Pineda is with the Department
of Electrical Engineering, University of M\'alaga, M\'alaga, Spain. E-mail: spinedamorente@gmail.com.}
\thanks{R. Fern\'andez-Blanco and J. M. Morales are with the Department of Applied Mathematics, University of M\'alaga, M\'alaga, Spain. E-mails: ricardo.fcarramolino@gmail.com; juan.morales@uma.es}
\thanks{This work was supported in part by the Spanish Ministry of Economy, Industry and Competitiveness through projects ENE2016-80638-R and ENE2017-83775-P; and in part by the European
Research Council (ERC) under the EU Horizon 2020 research and innovation
programme (grant agreement No. 755705) and the Research Program for Young Talented Researchers of the University of M\'{a}laga through project PPIT-UMA-B1-2017/18. The author thankfully acknowledges the computer resources, technical expertise and assistance provided by the SCBI (Supercomputing and Bioinformatics) center of the University of Malaga.}}

\maketitle

\begin{abstract}
The short-term  operation  of a power  system  is  usually planned  by  solving  a  day-ahead  unit  commitment  problem. Due to historical reasons, the commitment of the power generating units is decided over a time horizon typically consisting of the 24 hourly periods of a day. In this paper, we show that, as a result of the increasing  penetration of intermittent renewable generation, this somewhat arbitrary and artificial division of time may prove to be significantly suboptimal and counterproductive. Instead, we propose a time-adaptive day-ahead unit commitment formulation that better captures the net-demand variability throughout the day. The proposed formulation provides the commitment and dispatch of thermal generating units over a set of 24 time periods too, but with different duration. To do that, we use a clustering procedure  to  select  the duration of those adaptive time periods  taking into account the renewable generation and demand forecasts. Numerical results show that, without increasing the computational burden, the proposed time-adaptive unit  commitment  allows  for  a  more  efficient  use  of  the  system flexibility,  which  translates  into  a  lower  operating  cost  and  a higher penetration of renewable production than those achieved by  a  conventional  hourly  unit  commitment  problem.
\end{abstract}

\begin{IEEEkeywords}
Clustering techniques, economic dispatch, renewable generation, unit commitment.
\end{IEEEkeywords}

%
\IEEEpeerreviewmaketitle

\section*{Nomenclature}
The main symbols used throughout this paper are explained next. Others are defined as required.

\subsection{Indexes and sets}

\begin{ldescription}{$xxxxxx$}
\item[$g$] Conventional generating unit index.
\item[$t$, $\tau$] Time interval indices. 
\end{ldescription}

\subsection{Parameters}

\begin{ldescription}{$xxxxxxxxxx$}
\item[$C^{LS}$] Load shedding cost (\euro/MWh).
\item[$C^M_{gt}$ ] Marginal production cost of unit $g$ at time interval $t$ (\euro/MWh).
\item[$C_{gt}^{SU}$] Start-up cost of unit $g$ at time interval $t$ (\euro).
\item[$d_t$] Duration of time interval $t$ (h).
\item[$DT_g$ $(UT_g)$] Minimum down (up) time of  unit $g$ (h).
\item[$DT^0_g$ $(UT^0_g)$] Number of hours unit $g$ has been offline (online) prior to the first period (h).
\item[$DT_g^I$ $(UT_g^I)$] Number of hours unit $g$ must be initially offline (online) due to its minimum down (up) time constraint (h).
\item[$\widehat{DT}_{gt}$ $(\widehat{UT}_{gt})$] Dynamic minimum down (up) time of unit $g$ at time interval $t$ (\# time periods).
\item[$\widehat{DT}_g^E$ $(\widehat{UT}_g^E)$] Minimum down (up) time of unit $g$ at the end of the time horizon (\# time periods).
\item[$\widehat{DT}_g^I$ $(\widehat{UT}_g^I)$] Number of time periods unit $g$ must be initially offline (online) due to its minimum down (up) time constraint.
\item[$N_T$] Number of time periods.
\item[$\overline{P}^D_t$] Demand at time $t$ (MW).
\item[$\overline{P}^G_g$ $(\underline{P}^G_g)$] Maximum (minimum) production of thermal unit $g$ (MW).
\item[$\overline{P}^S$ ($\overline{P}^W$)] Installed capacity of solar (wind) generation (MW).
\item[$RD_g$ $(RU_g)$] Ramp-down (Ramp-up) limit of unit $g$ (MW/h).
\item[$SD_g$ $(SU_g)$] Shutdown (Start-up) ramp limit of unit $g$ (MW/h).
\item[$U^0_g$] Initial commitment state of unit $g$ (1 if online, and 0 otherwise).
\item[$\alpha^S$ $(\alpha^W)$] Percentage of yearly demand covered by solar (wind) generation  (\%).
\item[$\rho^S_{t}$ $(\rho^W_{t})$] Capacity factor of solar (wind) generation at time interval $t$ (p.u.).
\end{ldescription}

\subsection{Variables}

\begin{ldescription}{$xxxxxxx$}
\item[$p^G_{gt}$] Power output of unit $g$ at time interval $t$ (MW).
\item[$p^S_{t}$ $(p^W_{t})$] Power from solar (wind) generation at time interval $t$ (MW).
\item[$p^D_t$] Satisfied demand in time interval $t$ (MW).
\item[$s_{gt}^U$] Start-up cost of unit $g$ at time interval $t$ (\euro).
\item[$u_{gt}$ ] Binary variable that is equal to 1 if thermal unit $g$ is online at time interval $t$ and 0 otherwise.
\end{ldescription}

\section{Introduction}

\IEEEPARstart{T}{he} objective of the unit commitment (UC) problem is to determine the on/off status and production level of all generating units to satisfy the electricity demand at the minimum operating cost taking into account the system-wide technical constraints \cite{Sheble1994}. Although initially designed to centrally operate power systems, the UC problem is also widely used in deregulated environments to obtain the accepted bids and offers that maximize the social welfare while complying with technical constraints \cite{Richter2000}. The UC problem is commonly formulated as a mixed-integer quadratic optimization problem and a review of the main methods to solve it can be found in \cite{Padhy2004,Pandzic2013}. Solving a UC problem is computationally expensive because of its combinatorial nature and therefore, some authors have proposed methods to reduce its computational burden \cite{Carrion2006,Morales-Espana2013}. Due to the high penetration of fluctuating renewable energy, day-ahead decisions have to be made facing a significant level of uncertainty, which has led to stochastic formulations of the UC problem \cite{Zheng2014}.

Besides the more important role of uncertainty in the operation of power systems, the integration of renewable generation requires further modifications of the traditional unit commitment and economic dispatch tools \cite{Philbrick2011}. One of the potential changes currently under debate relates to the time resolution chosen to determine the day-ahead commitment and dispatch quantities of thermal generating units. In 2012, the Federal Energy Regulatory Commission (FERC) stated that ``hourly transmission
scheduling protocols (...) are insufficient to provide system operators with the flexibility to manage their system effectively and efficiently'' in the FERC Order 764 \cite{Ferc764}. For this reason, the FERC proposed 15-minute schedules. 

Following this line of argument, the authors in \cite{Pandzzic2014} investigate the impact of time resolution on the performance of the UC problem. They found out that 15-minute schedules lead to substantial savings through more efficient commitment and dispatch decisions at the expense of significantly increasing the computational needs. Similarly, reference \cite{Deane2014} solves the UC problem for time resolutions of 5, 15, 30, and 60 min. They concluded that, in spite of the increase in computational times, higher time resolutions show benefits over traditional hourly simulation in power systems with relatively few flexible resources. The authors in \cite{Kazemi2016} also analyze the effects of different time resolutions such as 60, 30, 15, 10, and 5 min on UC results. Their conclusion was that UC should be implemented in a higher time resolution to efficiently overcome the intra-hour variations of renewable generation. The UC results presented in \cite{Morales-Espana2017} show that the use of hourly intervals to model the production of thermal units is an approximation leading to costly operational decisions. Alternatively, the authors of \cite{Bakirtzis2014,Bakirtzis2017} propose a unified UC and economic dispatch modeling tool that considers a finer time resolution during the first hours of the scheduling horizon and a coarser time resolution during the last ones. Their approach was proven to provide adequate capacity and ramping capability to follow sudden changes of renewable generation. Unlike the aforementioned works,  
where the impact of finer resolutions is examined, authors in \cite{VomStein2017} propose to aggregate hourly time periods in order to reduce the computational complexity of the stochastic economic dispatch problem. 

In short, existing works on this topic concluded that using finer time resolutions lead to operating cost savings compared to the results from the conventional hourly UC problem (CH-UC.) However, these savings involve a significant increase in the computational burden. Note that increasing the number of time periods leads to an exponential raise of the computational time required to solve the UC problem to optimality. In fact, as stated in \cite{Sioshansi2008a}, ``ISOs cannot currently solve their commitment problems to complete optimality within the allotted timeframe.'' Therefore, although existing works have demonstrated the benefits of finer time resolutions to solve the UC problem, its implementation in current power systems is still unrealistic due to computational limitations. 

\begin{figure} \centering
\begin{tikzpicture}[scale=0.5]
	\begin{axis}[	
    width=\textwidth,
    height=8cm,
    xmin = 0,
    xmax = 24,
    ymin = 8000,
    ymax = 26000,
	legend style={at={(0.5,0.96)},anchor=north,legend cell align=left,legend columns=2},
	clip marker paths=true,	
	xlabel = Time (h),
	ylabel = Net demand (MW)]
	\addplot[line width=1pt,draw=black!50,dashed] table [x=t, y=NetDemand, col sep=comma] {caiso_220418.csv}; \addlegendentry{Data}
    \addplot[const plot,line width=1pt,draw=black] table [x=t, y=NetDemand_1h, col sep=comma] {caiso_220418.csv}; \addlegendentry{CH-UC}
	\addplot[domain=0:1,fill=gray!15, draw=none] {25900}\closedcycle;
	\addplot[domain=2:3,fill=gray!15, draw=none] {25900}\closedcycle;
	\addplot[domain=4:5,fill=gray!15, draw=none] {25900}\closedcycle;
	\addplot[domain=6:7,fill=gray!15, draw=none] {25900}\closedcycle;
	\addplot[domain=8:9,fill=gray!15, draw=none] {25900}\closedcycle;
	\addplot[domain=10:11,fill=gray!15, draw=none] {25900}\closedcycle;
	\addplot[domain=12:13,fill=gray!15, draw=none] {25900}\closedcycle;
	\addplot[domain=14:15,fill=gray!15, draw=none] {25900}\closedcycle;
	\addplot[domain=16:17,fill=gray!15, draw=none] {25900}\closedcycle;
	\addplot[domain=18:19,fill=gray!15, draw=none] {25900}\closedcycle;
	\addplot[domain=20:21,fill=gray!15, draw=none] {25900}\closedcycle;
	\addplot[domain=22:23,fill=gray!15, draw=none] {25900}\closedcycle;
	\addplot[line width=1pt,draw=black!50,dashed] table [x=t, y=NetDemand, col sep=comma] {caiso_220418.csv}; 
    \addplot[const plot,line width=1pt,draw=black] table [x=t, y=NetDemand_1h, col sep=comma] {caiso_220418.csv};
	\end{axis}	
\end{tikzpicture} \\[2mm]
\begin{tikzpicture}[scale=0.5]
	\begin{axis}[	
    width=\textwidth,
    height=8cm,
    xmin = 0,
    xmax = 24,
    ymin = 8000,
    ymax = 26000,
	legend style={at={(0.5,0.96)},anchor=north,legend cell align=left,legend columns=4},
	clip marker paths=true,	
	xlabel = Time (h),
	ylabel = Net demand (MW)]
	\addplot[line width=1pt,draw=black!50,dashed] table [x=t, y=NetDemand, col sep=comma] {caiso_220418.csv}; \addlegendentry{Data}
    \addplot[const plot,line width=1pt,draw=black] table [x=t, y=NetDemand_agg, col sep=comma] {caiso_220418.csv}; \addlegendentry{TA-UC}
	\addplot[domain=0:1,fill=gray!15, draw=none] {25900}\closedcycle;
	\addplot[domain=2:3,fill=gray!15, draw=none] {25900}\closedcycle;
	\addplot[domain=4:5,fill=gray!15, draw=none] {25900}\closedcycle;
	\addplot[domain=6:7,fill=gray!15, draw=none] {25900}\closedcycle;
	\addplot[domain=8:9,fill=gray!15, draw=none] {25900}\closedcycle;
	\addplot[domain=10:11,fill=gray!15, draw=none] {25900}\closedcycle;
	\addplot[domain=12:13,fill=gray!15, draw=none] {25900}\closedcycle;
	\addplot[domain=14:15,fill=gray!15, draw=none] {25900}\closedcycle;
	\addplot[domain=16:17,fill=gray!15, draw=none] {25900}\closedcycle;
	\addplot[domain=18:19,fill=gray!15, draw=none] {25900}\closedcycle;
	\addplot[domain=20:21,fill=gray!15, draw=none] {25900}\closedcycle;
	\addplot[domain=22:23,fill=gray!15, draw=none] {25900}\closedcycle;
	\addplot[line width=1pt,draw=black!50,dashed] table [x=t, y=NetDemand, col sep=comma] {caiso_220418.csv}; 
    \addplot[const plot,line width=1pt,draw=black] table [x=t, y=NetDemand_agg, col sep=comma] {caiso_220418.csv}; 
	\end{axis}	
\end{tikzpicture} \\[2mm]
\caption{Time aggregation for the net demand in CAISO on April 22$^{\text{nd}}$, 2018 with a temporal high-resolution (5 minutes). The upper and lower plot correspond to the time aggregation in low-resolution for the CH-UC and TA-UC, respectively.} \label{fig:agg}
\end{figure}
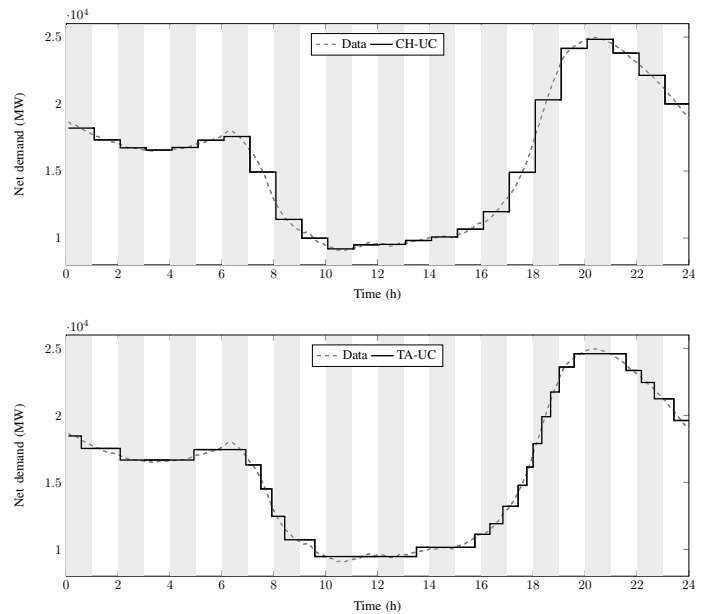

In order to overcome this drawback, this paper proposes a time-adaptive day-ahead UC formulation that makes a more efficient use of the system flexibility without increasing its computational burden. To illustrate the key idea of our proposal, Fig.~\ref{fig:agg} represents the net demand of CAISO \cite{CAISO} on April 22$^{\text{nd}}$, 2018 with a time resolution of 5 minutes, depicted in dashed lines in both plots.  More specifically, the upper plot of Fig. 1 represents the average net demand for each 1-hour time period in bold, as done in the CH-UC problem, thus leading to the conventional one-day planning horizon using 24 hourly time periods. The bold line in the lower plot of Fig.~1 represents the average net demand for 24 time periods of different duration according to the temporal aggregation proposed in this paper. We can notice that the way CH-UC approximates the real-time demand into 24 time periods is somewhat arbitrary and only used nowadays because of historical reasons. In fact, it is just apparent that significant errors are incurred by this approach to approximate the net demand variations for the time intervals 7:00--9:00 and 16:00--19:00.

Similarly to the conventional approach, our proposal also consists of approximating the one-day planning horizon using 24 time intervals. However, the duration of these 24 time periods is not arbitrarily set to 1 hour but determined following a rational criterion in order to better capture the dynamics of the time-dependent parameters such as the net demand. The bold line in the lower plot of Fig.~\ref{fig:agg} depicts the time aggregation proposed in this paper. It can be observed that the net demand is practically constant from 10:00 to 13:00 and therefore, only one time period with a duration of 3 hours is considered. On the other hand, the net demand grows from 11 GW at 16:00 to 22 GW at 19:00 and thus 9 time intervals are used to more accurately approximate this 3-hour time interval. In doing so, the proposed time-adaptive unit commitment (TA-UC) problem is expected to make a more efficient use of the available flexible resources thanks to a more accurate characterization of the net demand variations and, consequently, reduce the total operating costs.

The contributions of this paper are thus threefold:
\begin{itemize}
    \item[-] We introduce and formalize the notion of what we call \emph{time-adaptive unit commitment}, which considers a time horizon consisting of 24 time periods of different duration, instead of using the traditional 24 hourly time periods for solving the unit commitment problem. We propose to apply a hierarchical clustering algorithm to determine the 24 time periods with different duration that better capture the net system demand. 
    \item[-] We adapt the mixed-integer linear formulation for the day-ahead unit commitment problem to accommodate time periods of different duration. Specifically, the unit commitment formulation is made up of the same constraints but with substantial modifications on the technical limits of generating units such as start-up ramp rates, shutdown ramp rates, ramp-up and -down rates, as well as minimum up and down times.
    \item[-] Finally, we validate the results by showing substantial cost savings compared to the traditional unit commitment problem as well as a better utilization of the renewable production. We also show that the time-adaptive unit commitment formulation would be better suited in systems with a high share of solar production. More importantly, all those advantages come without increasing the computational burden. 
\end{itemize}

The rest of this paper is organized as follows. The clustering technique used to determine the duration of time periods is explained in Section \ref{sec:clustering}, while the proposed TA-UC problem is formulated in Section \ref{sec:formulation}. Section \ref{sec:comparison} presents the methodology used to compare the proposed approach with the conventional one. The results of an illustrative example and a more realistic case study are provided in Sections \ref{sec:example} and \ref{sec:casestudy}, respectively. Finally, conclusions are duly drawn in Section \ref{sec:conclusions}. 

\section{Time-period Aggregation}\label{sec:clustering}

The procedure to determine the duration of time periods for the proposed TA-UC problem is based on clustering techniques \cite{Hastie2009}. In particular, we use  hierarchical agglomerative clustering since its outcome is independent on the initialization of the algorithm and additional conditions on how the clusters are formed can be readily incorporated. In this paper we use an agglomerative hierarchical clustering based on Ward's method \cite{Ward1963}. This method recursively merges the two clusters that minimally increase the within-cluster variance. The methodology used in this paper and described below is based on the one proposed in \cite{Pineda2018}. 

To apply the clustering method to a short-term operational problem we need to define high- and low-resolution time periods as follows:
\begin{itemize}
    \item \textit{High-resolution} time periods are those with a duration of 5 or 10 minutes. The number of high-resolution time periods in one day, which is denoted by $N'$, is 288 or 144, respectively. In Fig. \ref{fig:agg}, the dashed lines represent the net demand in high-resolution.
    \item \textit{Low-resolution} time periods are referred to those time periods where the granularity is coarse. The number of low-resolution time periods in one day, which is denoted by $N$, is always equal to 24. In Fig. \ref{fig:agg}, the bold lines represent the approximation of the net demand using 24 low-resolution time periods.
\end{itemize}
Let $\mathbf{x}_i$ be a vector containing the normalized values in high-resolution of all time-dependent parameters at time point $i$ (with $i$ running from 1 to $N'$). For instance, $\mathbf{x}_i$ may have just one element corresponding to the aggregated net demand in the simplest case. The duration of the $N$ low-resolution time periods is determined as follows:

\begin{enumerate}
\item Set the initial number of clusters $n$ to the total number of \emph{high-resolution} time periods $N'$. That is, in this step, each data point $\mathbf{x}_i$ constitutes a cluster and $n = N'$.
\item Determine the centroid $\overline{\mathbf{x}}_I$ of each cluster $I$ as
\begin{equation}
\overline{\mathbf{x}}_I = \frac{1}{N_I} \sum_{i \in I} \mathbf{x}_i \label{eq:centroid}
\end{equation}
\item Compute the dissimilarity between each pair of adjacent clusters $I,J$ according to Ward's method using the equation \eqref{eq:distance}, where $N_I$ and $N_J$ are the number of elements in clusters $I$ and $J$, respectively.
\begin{equation}
D(I,J) = \frac{2N_IN_J}{N_I+N_J}  ||\overline{\mathbf{x}}_I-\overline{\mathbf{x}}_J||^2 \label{eq:distance}
\end{equation}
\item Merge the two closest \emph{adjacent} clusters $(I',J')$ according to the dissimilarity matrix, i.e., $(I',J') \in \text{argmin } D(I,J)$ subject to $J \in \mathcal{A}(I)$, where $\mathcal{A}(I)$ is the set of clusters adjacent to cluster $I$. Two clusters $I$ and $J$ are said to be adjacent if $I$ contains a high-resolution time period that is consecutive to a  high-resolution time period in $J$, or vice versa.
\item Update $n \leftarrow n - 1$.
\item If $n=N$, go to step 7). Otherwise go to step 2).
\item Determine the value of the parameters (e.g., the net demand) for each low-resolution time-period as its cluster's centroid $\overline{\mathbf{x}}_I$.
\item The number of \emph{high-resolution} time periods belonging to each final cluster determines the duration of each \emph{low-resolution} time period, which is denoted by $d_t$.
\end{enumerate}


\section{Time-adaptive Unit Commitment}\label{sec:formulation}
In this section we present the formulation of the proposed day-ahead TA-UC problem. For the sake of simplicity, this problem is formulated as a deterministic model. In other words, demand levels and renewable capacity factors are forecast values for the following day. Additionally, we assume that the forecasts are available in high-resolution (5 or 10 minutes.) Following common practice in the technical literature, the day-ahead UC problem is solved considering 24 low-resolution time periods.
For the traditional approach, the resolution of these time periods is equal to one hour regardless of the forecast net demand. Conversely, the low-resolution time periods for the proposed day-ahead UC differ from one day to another according to the forecast system conditions. Therefore, determining the duration of the low-resolution time periods according to the clustering procedure presented in Section \ref{sec:clustering} must be performed once a day right after the forecasts are issued.

We also formulate the UC assuming a single-bus system that disregards network constraints. Based on the model in \cite{Carrion2006}, the proposed TA-UC problem is formulated as the following mixed-integer linear program:

\begin{align}
& { \underset{\Xi}{\rm min}} \;  \sum_{g,t} \left( C_g^{M} p^G_{gt} d_t + s_{gt}^U \right) +  \sum_{t} C^{LS}d_t \left(\overline{P}^D_t - p^D_t\right)  \label{uc_of} \\
& {\text{subject to:}} \hspace{8cm} \nonumber \\
& \sum_g p^G_{gt} + p^W_{t} + p^S_{t} = p^D_t, \quad \forall t \label{uc_bal} \\
& 0 \leq p^D_t \leq \overline{P}^D_t, \quad \forall t \label{uc_dem} \\
& 0 \leq p^W_{t} \leq \rho^W_{t} \overline{P}^W, \quad \forall t \label{uc_wmax}\\
& 0 \leq p^S_{t} \leq \rho^S_{t} \overline{P}^S, \quad \forall t \label{uc_solarmax}\\
& u_{gt} \underline{P}^G_g \leq p^G_{gt} \leq u_{gt}\overline{P}^G_g, \quad \forall g, \forall t \label{uc_gmax} \\
%
%
& s_{gt}^U \geq C_g^{SU} \left( u_{gt}-u_{g,t-1} \right), \quad \forall g, \forall t \label{uc_SUcost1} \\
& s_{gt}^U \geq 0, \quad \forall g, \forall t \label{uc_SUcost2} \\
%
%
& p^G_{gt} - p^G_{g,t-1} \leq  \widehat{RU}_{gt} u_{g,t-1} +  \widehat{SU}_{gt} \left( u_{gt}-u_{g,t-1} \right)  \notag \\
& \hspace{2cm} + \overline{P}_g^G \left( 1-u_{gt}\right), \forall g, \forall t \label{uc_ramp1} \\
%
%
%
& p^G_{g,t-1} - p^G_{gt} \leq \widehat{RD}_{gt} u_{gt} +  \widehat{SD}_{gt} \left(u_{g,t-1}-u_{gt}\right)  \notag\\ 
& \hspace{2cm}  + \overline{P}_g^G \left(1-u_{g,t-1}\right), \forall g, \forall t \label{uc_ramp2}\\
& p^G_{gt} \leq \overline{P}^G_g u_{g,t+1} + \widehat{SD}_{gt} \left( u_{gt}-u_{g,t+1} \right), \forall g, \forall t < N_T \label{uc_ramp3} \\
& \sum_{t=1}^{\widehat{UT}^I_g} \left(1-u_{gt} \right)=0, \quad \forall g \label{uc_ut0}\\
& \sum_{\tau = t}^{t+\widehat{UT}_{gt}-1} u_{g\tau} \geq \widehat{UT}_{gt}\left( u_{gt} - u_{g,t-1} \right), \notag\\ 
& \hspace{2cm} \forall g, \forall t=\widehat{UT}^I_g + 1 \dots N_T - \widehat{UT}_g^E +1 \label{uc_ut1}\\
& \sum_{\tau=t}^{N_T} \left( u_{g\tau} - \left( u_{gt}-u_{g,t-1} \right) \right) \geq 0, \notag\\ 
& \hspace{2cm} \forall g, \forall t = N_T - \widehat{UT}_g^E +2 \dots N_T \label{uc_ut2} \\
& \sum_{t=1}^{\widehat{DT}^I_g} u_{gt} = 0, \quad \forall g \label{uc_dt0}\\
& \sum_{\tau=t}^{t+\widehat{DT}_{gt}-1} \left( 1-u_{g\tau}  \right) \geq \widehat{DT}_{gt} \left( u_{g,t-1}-u_{gt} \right), \notag\\ 
& \hspace{2cm} \forall g, \forall t = \widehat{DT}^I_g + 1 \dots N_T - \widehat{DT}_g^E +1 \label{uc_dt1} \\
& \sum_{\tau=t}^{N_T}\left( 1 - u_{g\tau} - \left( u_{g,t-1} - u_{gt} \right) \right) \geq 0, \notag\\ 
& \hspace{2cm} \forall g, \forall t = N_T - \widehat{DT}_g^E+2 \dots N_T \label{uc_dt2} \\
%
%
%
%
%
%
%
%
& u_{gt} \in \{0,1\},  \quad \forall g, \forall t, \label{uc_bin} 
\end{align}

\noindent where $\Xi$ is the set of  dispatch decisions ($p^G_{gt}$, $p^W_{t}$, $p^O_{t}$, $p^D_{t}$) and commitment decisions ($u_{gt}$,  $s^U_{gt}$).

The objective function \eqref{uc_of} comprises three terms, namely the production and start-up costs of thermal units as well as the costs due to load shedding. Constraints \eqref{uc_bal} model the power balance at each time interval. This balance includes injections from thermal and renewable generation. Constraints \eqref{uc_dem}--\eqref{uc_gmax} set bounds on satisfied demand, wind power production, solar power production, and thermal power production, respectively. Start-up costs are modeled by equations \eqref{uc_SUcost1}
and \eqref{uc_SUcost2}. Constraints \eqref{uc_ramp1}--\eqref{uc_ramp3} enforce the ramping rates of thermal units as described in \cite{Carrion2006}, where the ramping rates in MW are given as follows:
\begin{align}
    & \widehat{RU}_{gt} = \min\{\overline{P}^G_g, \max\{\underline{P}^G_g, RU_g \widehat{d}_t\}\}, \quad \forall g, \forall t \\
    & \widehat{RD}_{gt} = \min\{\overline{P}^G_g, \max\{\underline{P}^G_g, RD_g \widehat{d}_t\}\}, \quad \forall g, \forall t \\
    & \widehat{SU}_{gt} = \min\{\overline{P}^G_g, \max\{\underline{P}^G_g, SU_g \widehat{d}_t\}\}, \quad \forall g, \forall t \\
    & \widehat{SD}_{gt} = \min\{\overline{P}^G_g, \max\{\underline{P}^G_g, SD_g \widehat{d}_t\}\}, \quad \forall g, \forall t, 
\end{align}

\noindent where $\widehat{d}_t = 0.5(d_{t-1} + d_t)$, i.e., the power output is assumed to ramp up or down from the middle point of intervals $t-1$ and $t$. Note that if the conventional approach is used, all low-resolution time periods are one-hour long and the ramping rates are the same for all time periods. On the other hand, if the low-resolution time periods have different duration, the ramping rate corresponding to a short-time period is lower than that corresponding to a long-time period. In other words, the longer the time period duration, the larger the allowed variation of the power output.

Constraints \eqref{uc_ut0}--\eqref{uc_ut2} and \eqref{uc_dt0}--\eqref{uc_dt2} correspond to the minimum up-time and down-time constraints of thermal units, in that order. These constraints allow for predefined time intervals of different duration, by computing parameters $\widehat{UT}_g^I$, $\widehat{UT}_g^E$, $\widehat{UT}_{gt}$, $\widehat{DT}_g^I$, $\widehat{DT}_g^E$, and $\widehat{DT}_{gt}$ \emph{ex-ante} as follows:
\begin{align}
& \widehat{UT}_g^I =  \underset{\omega=1,2,\ldots}{\arg\min} \; \omega : \sum_{\tau=1}^{\omega}d_{\tau} \geq UT_g^I, \quad \forall g \\
& \widehat{UT}_g^E =  \underset{\omega=1,2,\ldots}{\arg\min} \; \omega : \sum_{\tau=N_T-\omega+1}^{N_T}d_{\tau} \geq UT_g, \quad \forall g \\
& \widehat{UT}_{gt} =  \underset{\omega=1,2,\ldots}{\arg\min} \; \omega : \sum_{\tau=t}^{t+\omega-1}d_{\tau} \geq UT_g, \quad \forall g, \forall t\\
& \widehat{DT}_g^I =  \underset{\omega=1,2,\ldots}{\arg\min} \; \omega : \sum_{\tau=1}^{\omega}d_{\tau} \geq DT_g^I, \quad \forall g \\
& \widehat{DT}_g^E =  \underset{\omega=1,2,\ldots}{\arg\min} \; \omega : \sum_{\tau=N_T-\omega+1}^{N_T}d_{\tau} \geq DT_g, \quad \forall g \\
& \widehat{DT}_{gt} =  \underset{\omega=1,2,\ldots}{\arg\min} \; \omega : \sum_{\tau=t}^{t+\omega-1}d_{\tau} \geq DT_g, \quad \forall g, \forall t
\end{align}
\noindent where $UT_g^I = \min \{ N_T, \left( UT_g - UT_g^0 \right) U_g^0 \}$ and $DT_g^I = \min \{ N_T, \left( DT_g - DT_g^0 \right) U_g^0 \}$. Note that for the conventional UC, the duration of all low-resolution time periods is equal to one hour, i.e.,  $d_t=1$, $\forall t$. Therefore, the minimum up- and down times are the same for all time periods of the following day, i.e., $\widehat{UT}_{gt}=UT_g$, $\forall t$ and $\widehat{DT}_{gt}=DT_g$, $\forall t$. If the duration of time periods changes throughout the day according to the proposed approach, then the minimum number of time periods that a generating unit has to be on or off may also change.
%
Finally, the integrality of binary variables is imposed by constraints \eqref{uc_bin}. 

To sum up, formulation \eqref{uc_of}--\eqref{uc_bin} is valid both for the conventional UC problem, in which the duration of the 24 low-resolution time periods is equal to one hour, and for the proposed time-adaptive UC problem, in which the duration of the time periods is determined according to the clustering procedure presented in Section \ref{sec:clustering}. In the former case, the technical limits of the thermal generating units such as ramping limits or minimum up and down times are the same for all time periods and optimization model \eqref{uc_of}--\eqref{uc_bin} becomes the CH-UC in  \cite{Carrion2006}. In the latter case, such technical limits are different for each time period depending on their duration and model \eqref{uc_of}--\eqref{uc_bin} becomes the proposed TA-UC.

For the sake of simplicity, the UC problem solved in this paper disregards network constraints. Therefore, consecutive time periods are merged based on the aggregated net demand for the whole system. If a transmission network is included, the distance between two consecutive time periods can be computed by comparing the vectors containing the net demand at each network bus by using, for example, the Euclidean norm. In this way, two consecutive time periods would be similar if the net demand across all network buses is close enough.

Note that we also assume a centralized or cost-based environment rather than a market-based system. We do this in order to focus the discussion on the main contribution of this paper, which is the proposal of a UC model that considers 24 time periods of different duration. Finally, it is worth mentioning that formulation \eqref{uc_of}--\eqref{uc_bin} preclude balancing mechanisms such as reserves in order to assess the performance of the proposed time-period aggregation approach using a simple modeling of reality.

\section{Comparison Methodology} \label{sec:comparison}

In this section, we present the procedure to compare CH-UC and TA-UC. First, it is worth clarifying that thermal generating units are divided into three main groups according to their flexibility: base-load units, medium-load units, and peak-load units. This implies that the commitment and dispatch of the base-load units and the commitment of the medium-load units have to be necessarily decided one day in advance and cannot be modified in the real-time operation of the system. On the other hand, the dispatch of the medium-load units and the commitment and dispatch of the peak-load units can adapt to the real-time net-demand. That said, the total operating cost for each methodology is computed as follows:

\begin{enumerate}
    \item Forecast values for time-dependent parameters such as demand level and renewable capacity factors are provided in a high-resolution time scale (5 or 10 min.)
    \item Original data are approximated using 24 low-resolution time periods whose duration depends on each approach:
    \begin{enumerate}
        \item CH-UC: the duration of 
        all low-resolution time periods is equal to one hour.
        \item TA-UC: the duration of the low-resolution time periods is determined as explained in Section \ref{sec:clustering}.
    \end{enumerate}
    \item Model \eqref{uc_of}--\eqref{uc_bin} is solved for the low-resolution time series computed as described in step 2). This problem provides the day-ahead commitment and dispatch of base-load units and the day-ahead commitment of medium-load units.
    \item Model \eqref{uc_of}--\eqref{uc_bin} is solved again considering the original high-resolution time data and with the commitment and dispatch of base-load units and the commitment of medium-load units fixed to those obtained in step 3.) This way, we are able to simulate the real-time operation of the system and the performance of the day-ahead decisions given by each methodology.
\end{enumerate}

Let us denote the daily operating cost obtained in step 4) for CH-UC and TA-UC as $C^{CH}$ and $C^{TA}$, respectively. We evaluate the performance of the proposed approach by computing the relative difference between $C^{CH}$ and $C^{TA}$ as:
\begin{equation}
    \Delta C (\%) = 100\cdot\frac{C^{CH} - C^{TA}}{C^{CH}}
\end{equation}

In order to draw conclusions about the proposed methodology, the daily cost difference $\Delta C$ must be computed for several consecutive days. To this end, the UC model \eqref{uc_of}--\eqref{uc_bin} is run in a rolling horizon with an 8-hour look-ahead window. The dispatch and commitment decisions in day $D$ are obtained by running the model in day $D$ plus a look-ahead window of length equal to 8 time intervals, which corresponds to the next day $D+1$. The initial conditions are taken from the previous high-resolution simulation at the end of the time span from day $D-1$. Finally, it is worth mentioning that the two approaches compared in this paper involve computational burdens of the same order of magnitude since both include the same number of constraints and continuous and binary variables.

\section{Illustrative Example} \label{sec:example}

This section illustrates the performance of the proposed TA-UC by using a stylized example of six generating units, whose data is collated in Table \ref{tab:gen_exa}. For the sake of simplicity, minimum up and down times and ramp limits are neglected. Besides, it is assumed that production and start-up costs of thermal units remain unchanged over the time horizon.

\renewcommand{\arraystretch}{1.2}
\begin{table}
    \centering
    \caption{\textsc{Generating Unit Data -- Example}}
    \begin{tabular}{ccccc}
    \hline
    Technology & $\underline{P}_g^G$ (MW) & $\overline{P}_g^G$ (MW) & $C_{gt}^M$ (\euro/MWh) & \# units \\
    \hline
    Base & 150 & 200 & 10 & 4 \\
    Medium & 50 & 100 & 30 & 1 \\
    Peak & 0 & 50 & 50 & 1 \\
    \hline
    \end{tabular} 
    \label{tab:gen_exa}
\end{table}

The time horizon of this illustrative example spans six time periods of 30 minutes each. Table \ref{tab:dem_exam} provides the demand, the solar power production as well as the net demand for each time period. Observe that the demand and the solar production increases and decreases in the last two time periods, respectively. The load shedding cost is set to \euro100/MWh.

\begin{table}
    \centering
    \caption{\textsc{Net Demand -- Example}}
    \begin{tabular}{ccccccc}
    \hline
    Time period & $t1$ & $t2$ & $t3$ & $t4$ & $t5$ & $t6$ \\
    \hline
    Duration (h) & 0.5 & 0.5 & 0.5 & 0.5 & 0.5 & 0.5 \\
    Demand (MW) & 500 & 500 & 500 & 500 & 650 & 850 \\
    Solar (MW) & 300 & 300 & 300 & 300 & 200 & 0 \\
    Net demand (MW) & 200 & 200 & 200 & 200 & 450 & 850 \\
    \hline
    \end{tabular}
    \label{tab:dem_exam}
\end{table}

To determine the commitment of units in a conventional fashion, the 30-minute time periods are merged two by two in order to determine the hourly day-ahead commitment and dispatch of each unit, provided in the upper part of Table \ref{tab:cuc_exam}. The real-time generation levels of each generating technology, the load shedding, and the solar spillage for each 30-minute time period of the optimization horizon are shown in the lower part of the same table. Although the net demand of $t5$ and $t6$ are quite different, the day-ahead decisions for these two time periods are the same under this approach. This involves some solar spillage in $t5$, the start-up of the expensive peak unit and some load shedding in $t6$. The total operating cost for this conventional UC plan is 18500\euro. 

\begin{table}
    \centering
    \caption{\textsc{Conventional Hourly Unit Commitment (CH-UC) --  Example}}
    \begin{tabular}{ccccccc}
    \hline
    \multicolumn{7}{c}{Day-ahead dispatch} \\
    \hline
    Time periods & \multicolumn{2}{c}{$t1+t2$} & \multicolumn{2}{c}{$t3+t4$} & \multicolumn{2}{c}{$t5+t6$} \\
    \hline
    Net demand (MW) & \multicolumn{2}{c}{200} & \multicolumn{2}{c}{200} &  \multicolumn{2}{c}{650} \\
    Base (MW) & \multicolumn{2}{c}{200} & \multicolumn{2}{c}{200} &  \multicolumn{2}{c}{600} \\
    Medium (MW) & \multicolumn{2}{c}{0} & \multicolumn{2}{c}{0} &  \multicolumn{2}{c}{50} \\
    Peak (MW) & \multicolumn{2}{c}{0} & \multicolumn{2}{c}{0} &  \multicolumn{2}{c}{0} \\
    \hline 
    \hline
    \multicolumn{7}{c}{Real-time operation} \\
    \hline
    Time periods & $t1$ & $t2$ & $t3$ & $t4$ & $t5$ & $t6$ \\
    \hline
    Net demand (MW) & 200 & 200 & 200 & 200 & 450 & 850 \\
    Base (MW) & 200 & 200 & 200 & 200 & 600 & 600 \\
    Medium (MW) & 0 & 0 & 0 & 0 & 50 & 100 \\
    Peak (MW) & 0 & 0 & 0 & 0 & 0 & 50 \\
    Load shed (MW) & 0 & 0 & 0 & 0 & 0 & 100 \\
    Solar spillage (MW) & 0 & 0 & 0 & 0 & 200 & 0 \\
    \hline
    \end{tabular}
    \label{tab:cuc_exam}
\end{table}

The TA-UC proposed in this paper also considers three time periods to determine the day-ahead commitment and dispatch of generating units. However, the time period aggregation is quite different. Since the first four 30-minute time periods have the same net demand, they are merged into a 2-hour time period, as illustrated in the upper part of Table \ref{tab:auc_exam}. In this way, the day-ahead decisions adapt better to the net demand changes happening in $t5$ and $t6$ and load shedding and solar spillage are no longer required in the real-time operation. Under this approach, the total cost amounts to 11500\euro, then involving a 38\% cost reduction with respect to the CH-UC. Due to its simplicity, the computational times of this illustrative example for both methods are negligible.
   
\begin{table}
    \centering
    \caption{\textsc{Time-Adaptive Unit Commitment (TA-UC) -- Example}}
    \begin{tabular}{ccccccc}
    \hline
    \multicolumn{7}{c}{Day-ahead dispatch} \\
    \hline
    Time periods & \multicolumn{4}{c}{$t1+t2+t3+t4$} & \multicolumn{1}{c}{$t5$} & \multicolumn{1}{c}{$t6$} \\
    \hline
    Net demand (MW) & \multicolumn{4}{c}{200} & \multicolumn{1}{c}{450} &  \multicolumn{1}{c}{850} \\
    Base (MW) & \multicolumn{4}{c}{200} & \multicolumn{1}{c}{400} &  \multicolumn{1}{c}{800} \\
    Medium (MW) & \multicolumn{4}{c}{0} & \multicolumn{1}{c}{50} &  \multicolumn{1}{c}{50} \\
    Peak (MW) & \multicolumn{4}{c}{0} & \multicolumn{1}{c}{0} &  \multicolumn{1}{c}{0} \\
    \hline 
    \hline
    \multicolumn{7}{c}{Real-time operation} \\
    \hline
    Time periods & $t1$ & $t2$ & $t3$ & $t4$ & $t5$ & $t6$ \\
    \hline
    Net demand (MW) & 200 & 200 & 200 & 200 & 450 & 850 \\
    Base (MW) & 200 & 200 & 200 & 200 & 400 & 800 \\
    Medium (MW) & 0 & 0 & 0 & 0 & 50 & 50 \\
    Peak (MW) & 0 & 0 & 0 & 0 & 0 & 0 \\
    Load shed (MW) & 0 & 0 & 0 & 0 & 0 & 0 \\
    Solar spillage (MW) & 0 & 0 & 0 & 0 & 0 & 0 \\
    \hline
    \end{tabular}
    \label{tab:auc_exam}
\end{table} 

\section{Case Study} \label{sec:casestudy}

The proposed methodology is tested using a more realistic case study based on the Spanish power system. The Spanish electricity demand (after subtracting the hydro power production) during the year 2017 in a 10-minute time resolution is scaled down a factor of 10 in order to keep the computational burden of the UC problems within reasonable limits. The capacity factors of solar and wind power production are also taken from the Spanish system during the year 2017 in a 10-minute resolution. These data are publicly available in \cite{esios}. The installed capacity of wind and solar power generation are determined by equations \eqref{eq:installed} so that the share of demand covered by these technologies amounts to $\alpha^S$ and $\alpha^W$, respectively.
\begin{equation}
    \overline{P}^S = \alpha^S \frac{\sum_t \overline{P}^D_t}{ \sum_t \rho^S_t} \qquad \qquad \overline{P}^W = \alpha^W \frac{\sum_t \overline{P}^D_t}{ \sum_t \rho^W_t} \label{eq:installed}
\end{equation}

The generation portfolio is composed of 3 base-load units, 4 medium-load units, and 6 peak-load units, whose technical and economic data are provided in Table \ref{tab:gen_exa2}. These parameters have been chosen according to \cite{DeJonghe2014a} and the references therein. As observed, base-load units are cheap albeit inflexible, while peak-load power plants are flexible albeit expensive units. Note that the ramp rates, minimum up and down times, and marginal costs of units of the same type are slightly different in order to represent the variety of generating technologies. In all cases, higher flexibility implies higher marginal costs. Besides, it is assumed that production and start-up costs of thermal units remain unchanged over the time horizon. Moreover, we assume that the minimum power output of peak-load units is 0 and their minimum up and down times are disregarded. Finally, the cost of load shedding is set to \euro10000/MWh.

\renewcommand{\arraystretch}{1.2}
\setlength{\tabcolsep}{1.5mm}
\begin{table}
    \caption{\textsc{Thermal Generating Unit Data -- Case Study}}
    \centering
    \begin{tabular}{cccccccccc}
    \hline
Unit & Type & $\underline{P}_g^G$ & $\overline{P}_g^G$ & $RU_g$ & $RD_g$ & $UT_g$ & $DT_g$ & $C^{SU}_{gt}$ & $C^M_{gt}$ \\
    \hline
    $g1$ & Base & 200 & 400 & 120 & 120 & 9 & 9 & 1000 & 20 \\
    $g2$ & Base & 200 & 400 & 130 & 130 & 8.5 & 8.5 & 1000 & 21 \\
    $g3$ & Base & 200 & 400 & 140 & 140 & 8 & 8 & 1000 & 22 \\
    $g4$ & Medium & 100 & 300 & 105 & 105 & 5 & 5 & 800 & 50 \\
    $g5$ & Medium & 100 & 300 & 120 & 120 & 4.7 & 4.7 & 800 & 51 \\
    $g6$ & Medium & 100 & 300 & 135 & 135 & 4.3 & 4.3 & 800 & 52 \\
    $g7$ & Medium & 100 & 300 & 150 & 150 & 4 & 4 & 800 & 53 \\
    $g8$ & Peak & 0 & 250 & 125 & 125 & - & - & 500 & 80 \\
    $g9$ & Peak & 0 & 250 & 130 & 130 & - & - & 500 & 81 \\
    $g10$ & Peak & 0 & 250 & 135 & 135 & - & - & 500 & 82 \\
    $g11$ & Peak & 0 & 250 & 140 & 140 & - & - & 500 & 83 \\
    $g12$ & Peak & 0 & 250 & 145 & 145 & - & - & 500 & 84 \\
    $g13$ & Peak & 0 & 250 & 150 & 150 & - & - & 500 & 85 \\
    \hline
    \end{tabular} 
    \label{tab:gen_exa2}
\end{table}

The proposed TA-UC and conventional CH-UC models have been implemented on a Linux-based server with one CPU clocking at 2.6 GHz and 20 GB of RAM using CPLEX 12.6.3 \cite{Cplex} under Pyomo 5.2 \cite{Pyomo}. Optimality gap is set to 0\%. 

For the sake of illustration, first we present  the day-ahead schedule provided by each UC formulation for two typical days of 2017, namely March 19$^\text{th}$ and December 26$^\text{th}$. The yearly penetration of wind and solar generation is 20\% and 20\%, respectively. March 19$^\text{th}$ is thus a representative day with abrupt intra-day net demand variations, which gives rise to the so-called ``duck curve,'' whereas December 26$^\text{th}$ is characterized by smooth intra-day net demand variations.   

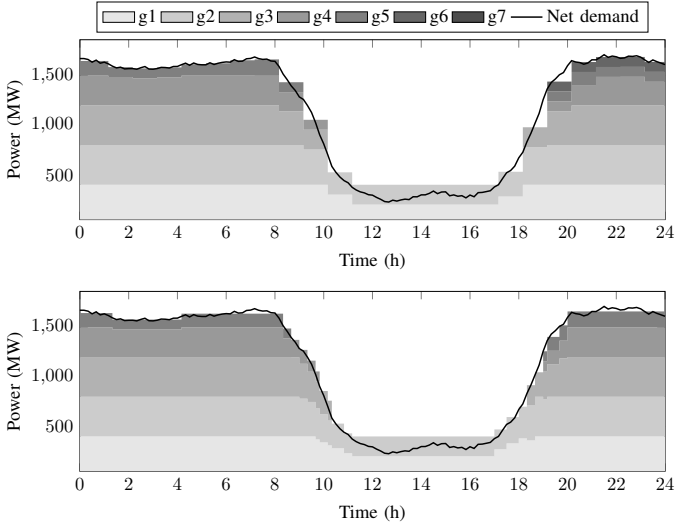
\begin{figure} \centering
\begin{tikzpicture}[scale=0.7]
	\begin{axis}[	
    width=0.7\textwidth,
    height=5cm,
    xmin = 0,
    xmax = 24,
	legend style={at={(0.5,1.04)},anchor=south,legend cell align=left,legend columns=8},
	xlabel = Time (h),
	ylabel = Power (MW)
    ]		
	\addplot[draw=none, const plot,stack plots=y,area style,enlarge x limits=false,fill=black!10] table [x=t, y=g1, col sep=comma] {78hc.csv}\closedcycle;\addlegendentry{g1}
	\addplot[draw=none,const plot,stack plots=y,area style,enlarge x limits=false,fill=black!20] table [x=t, y=g2, col sep=comma] {78hc.csv}
	\closedcycle; \addlegendentry{g2}
	\addplot[draw=none,const plot,stack plots=y,area style,enlarge x limits=false,fill=black!30] table [x=t, y=g3, col sep=comma] {78hc.csv}
	\closedcycle; \addlegendentry{g3}
	\addplot[draw=none,const plot,stack plots=y,area style,enlarge x limits=false,fill=black!40] table [x=t, y=g4, col sep=comma] {78hc.csv}
	\closedcycle; \addlegendentry{g4}
	\addplot[draw=none,const plot,stack plots=y,area style,enlarge x limits=false,fill=black!50] table [x=t, y=g5, col sep=comma] {78hc.csv}
	\closedcycle; \addlegendentry{g5}
	\addplot[draw=none,const plot,stack plots=y,area style,enlarge x limits=false,fill=black!60] table [x=t, y=g6, col sep=comma] {78hc.csv}
	\closedcycle; \addlegendentry{g6}
	\addplot[draw=none,const plot,stack plots=y,area style,enlarge x limits=false,fill=black!70] table [x=t, y=g7, col sep=comma] {78hc.csv}
	\closedcycle; \addlegendentry{g7}
	\addplot[line width=0.8pt,draw=black] table [x=t, y=d, col sep=comma] {78hc.csv}; \addlegendentry{Net demand}
	\end{axis}	
\end{tikzpicture} \\[2mm]
\begin{tikzpicture}[scale=0.7]
	\begin{axis}[	
    width=0.7\textwidth,
    height=5cm,
    xmin = 0,
    xmax = 24,
	xlabel = Time (h),
	ylabel = Power (MW)
    ]		
	\addplot[draw=none, const plot,stack plots=y,area style,enlarge x limits=false,fill=black!10] table [x=t, y=g1, col sep=comma] {78ta.csv}\closedcycle;
	\addplot[draw=none,const plot,stack plots=y,area style,enlarge x limits=false,fill=black!20] table [x=t, y=g2, col sep=comma] {78ta.csv}
	\closedcycle; 
	\addplot[draw=none,const plot,stack plots=y,area style,enlarge x limits=false,fill=black!30] table [x=t, y=g3, col sep=comma] {78ta.csv}
	\closedcycle; 
	\addplot[draw=none,const plot,stack plots=y,area style,enlarge x limits=false,fill=black!40] table [x=t, y=g4, col sep=comma] {78ta.csv}
	\closedcycle; 
	\addplot[draw=none,const plot,stack plots=y,area style,enlarge x limits=false,fill=black!50] table [x=t, y=g5, col sep=comma] {78ta.csv}
	\closedcycle; 
	\addplot[draw=none,const plot,stack plots=y,area style,enlarge x limits=false,fill=black!60] table [x=t, y=g6, col sep=comma] {78ta.csv}
	\closedcycle; 
	\addplot[draw=none,const plot,stack plots=y,area style,enlarge x limits=false,fill=black!70] table [x=t, y=g7, col sep=comma] {78ta.csv}
	\closedcycle;
	\addplot[line width=0.8pt,draw=black] table [x=t, y=d, col sep=comma] {78ta.csv};
	\end{axis}	
\end{tikzpicture} 
\caption{Day-ahead dispatch and real-time net demand for the CH-UC (upper plot) and the TA-UC (lower plot) on March 19$^\text{th}$, 2017.} \label{fig:march19}
\end{figure}

The bold lines in Fig. \ref{fig:march19} plot the real-time net demand on March 19$^\text{th}$, 2017 in high-resolution. The stacked area chart represent the day-ahead dispatch of base-load units (g1, g2 and g3) and medium-load units (g4, g5, g6 and g7) in low-resolution. The upper subplot corresponds to the results from CH-UC and therefore the duration of all low-resolution time periods is equal to one hour. The lower subplot to the ones from TA-UC and as observed, the duration of the low-resolution time periods varies throughout the day. If the stacked chart exceeds the bold line, the day-ahead dispatch of base- and medium-load units is higher than the net demand and then, wind or solar spillage will occur in real-time operation. Conversely, if the bold curve surpasses the stacked chart, the net demand is higher than the day-ahead dispatch of base- and medium-load units. Such mismatch has to be covered by (i) increasing the power output of online medium-load units whenever possible, (ii) using more expensive peak-load units, or (iii) resorting to load shedding. In either case, the larger the differences between day-ahead dispatch and real-time net demand, the more expensive the real-time operation of the system becomes. Although not considered in this work, the differences between generation and consumption can be partially compensated with other expensive real-time balancing resources such as frequency control.

We can observe that the net demand of Fig. \ref{fig:march19} has two flat peaks (from 0:00 to 8:00 and from 20:00 to 23:59), a flat valley (from 11:00 to 17:00), and very steep net demand variations between the peaks and the valley (from 8:00 to 11:00 and from 17:00 to 20:00). This net demand profile is caused by a high solar production and a low wind penetration. The proposed model TA-UC takes advantage of that and chooses very long time periods during the peaks and the valley, but shorter time periods in the transitions, as can be seen in the lower subplot. Comparing the two subplots of Fig. \ref{fig:march19}, it can be observed that the proposed TA-UC approximates the net demand much better than the conventional approach, which translates into a lower real-time operating cost. In fact, the real-time cost for TA-UC and CH-UC is \euro765366 and \euro783643, respectively. The relative cost saving for this day  amounts to 2.33\%. In practical terms, the lower mismatch between generation and consumption would entail a reduction of expensive frequency control actions.

Table \ref{tab:march19} provides the daily shares of the different production types on March 19$^\text{th}$, 2017. It can be noticed that the TA-UC achieves  higher shares from renewable and base generation (87.87\%) than the CH-UC (87.26\%). This implies a lower utilization of the more expensive and polluting medium- and peak-load generating units.

\renewcommand{\arraystretch}{1.2}
\setlength{\tabcolsep}{1.5mm}
\begin{table}
    \centering
    \caption{\textsc{Daily Generation Shares on March} 19$^\text{th}$, 2017}
    \begin{tabular}{cccccc}
    \hline
Model & Wind (\%) & Solar (\%) & Base (\%) & Medium (\%) & Peak (\%) \\
    \hline
    CH-UC & 9.52 & 32.48 & 45.26 & 12.36 & 0.38 \\
    TA-UC & 9.47 & 32.87 & 45.53 & 12.06 & 0.07 \\
    \hline
    \end{tabular} 
    \label{tab:march19}
\end{table}

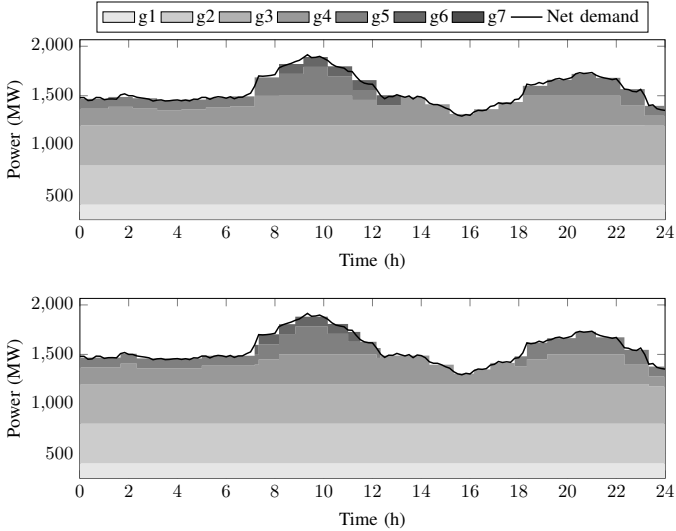
\begin{figure} \centering
\begin{tikzpicture}[scale=0.7]
	\begin{axis}[	
    width=0.7\textwidth,
    height=5cm,
    xmin = 0,
    xmax = 24,
	legend style={at={(0.5,1.04)},anchor=south,legend cell align=left,legend columns=8},	
	xlabel = Time (h),
	ylabel = Power (MW)
    ]		
	\addplot[draw=none, const plot,stack plots=y,area style,enlarge x limits=false,fill=black!10] table [x=t, y=g1, col sep=comma] {360hc.csv}\closedcycle;\addlegendentry{g1}
	\addplot[draw=none,const plot,stack plots=y,area style,enlarge x limits=false,fill=black!20] table [x=t, y=g2, col sep=comma] {360hc.csv}
	\closedcycle; \addlegendentry{g2}
	\addplot[draw=none,const plot,stack plots=y,area style,enlarge x limits=false,fill=black!30] table [x=t, y=g3, col sep=comma] {360hc.csv}
	\closedcycle; \addlegendentry{g3}
	\addplot[draw=none,const plot,stack plots=y,area style,enlarge x limits=false,fill=black!40] table [x=t, y=g4, col sep=comma] {360hc.csv}
	\closedcycle; \addlegendentry{g4}
	\addplot[draw=none,const plot,stack plots=y,area style,enlarge x limits=false,fill=black!50] table [x=t, y=g5, col sep=comma] {360hc.csv}
	\closedcycle; \addlegendentry{g5}
	\addplot[draw=none,const plot,stack plots=y,area style,enlarge x limits=false,fill=black!60] table [x=t, y=g6, col sep=comma] {360hc.csv}
	\closedcycle; \addlegendentry{g6}
	\addplot[draw=none,const plot,stack plots=y,area style,enlarge x limits=false,fill=black!70] table [x=t, y=g7, col sep=comma] {360hc.csv}
	\closedcycle; \addlegendentry{g7}
	\addplot[line width=0.8pt,draw=black] table [x=t, y=d, col sep=comma] {360hc.csv}; \addlegendentry{Net demand}
	\end{axis}	
\end{tikzpicture} \\[2mm]
\begin{tikzpicture}[scale=0.7]
	\begin{axis}[	
    width=0.7\textwidth,
    height=5cm,
    xmin = 0,
    xmax = 24,
	xlabel = Time (h),
	ylabel = Power (MW)
    ]		
	\addplot[draw=none, const plot,stack plots=y,area style,enlarge x limits=false,fill=black!10] table [x=t, y=g1, col sep=comma] {360ta.csv}\closedcycle;
	\addplot[draw=none,const plot,stack plots=y,area style,enlarge x limits=false,fill=black!20] table [x=t, y=g2, col sep=comma] {360ta.csv}
	\closedcycle; 
	\addplot[draw=none,const plot,stack plots=y,area style,enlarge x limits=false,fill=black!30] table [x=t, y=g3, col sep=comma] {360ta.csv}
	\closedcycle; 
	\addplot[draw=none,const plot,stack plots=y,area style,enlarge x limits=false,fill=black!40] table [x=t, y=g4, col sep=comma] {360ta.csv}
	\closedcycle; 
	\addplot[draw=none,const plot,stack plots=y,area style,enlarge x limits=false,fill=black!50] table [x=t, y=g5, col sep=comma] {360ta.csv}
	\closedcycle; 
	\addplot[draw=none,const plot,stack plots=y,area style,enlarge x limits=false,fill=black!60] table [x=t, y=g6, col sep=comma] {360ta.csv}
	\closedcycle; 
	\addplot[draw=none,const plot,stack plots=y,area style,enlarge x limits=false,fill=black!70] table [x=t, y=g7, col sep=comma] {360ta.csv}
	\closedcycle;
	\addplot[line width=0.8pt,draw=black] table [x=t, y=d, col sep=comma] {360ta.csv};
	\end{axis}	
\end{tikzpicture} 
\caption{Day-ahead dispatch and real-time net demand for the CH-UC (upper plot) and the TA-UC (lower plot) on December 26$^\text{th}$, 2017.} \label{fig:december26}
\end{figure}

Fig. \ref{fig:december26} plots the same results for December 26$^\text{th}$, 2017, which is a day with a particularly high wind production but barely any solar production. For this reason, the intra-day net demand variations are much smoother and the day-ahead dispatch of the TA-UC (lower subplot) approximates the net demand very similarly to the CH-UC (upper subplot). In fact, the real-time operating cost  for TA-UC and CH-UC is \euro1037129 and \euro1037257, which are practically the same. Table \ref{tab:december26} provides the daily shares of the different production types on December 26$^\text{th}$, 2017. As expected, all shares are basically the same for both models. Besides, the absence of very steep net demand variations makes the use of peak-load units unnecessary.

\renewcommand{\arraystretch}{1.2}
\setlength{\tabcolsep}{1.5mm}
\begin{table}
    \centering
    \caption{\textsc{Daily Generation Shares on December} 26$^\text{th}$, 2017}
    \begin{tabular}{cccccc}
    \hline
Model & Wind (\%) & Solar (\%) & Base (\%) & Medium (\%) & Peak (\%) \\
    \hline
    CH-UC & 38.35 & 3.30 & 44.99 & 13.35 & 0.00 \\
    TA-UC & 38.37 & 3.30 & 44.96 & 13.36 & 0.00 \\
    \hline
    \end{tabular} 
    \label{tab:december26}
\end{table}

Still maintaining the yearly penetration of wind and solar generation to 20\% each, we compare the annual results from both models. We remove January 1$^\text{st}$ and December 31$^\text{st}$ to neglect the impact of boundary conditions. Under these circumstances, the proposed TA-UC yields a yearly cost saving of 0.50\% compared to the results from the conventional CH-UC model. The proposed TA-UC outperforms the conventional CH-UC during 332 days, both models achieve the same operating cost in 9 days, and CH-UC provides lower operating costs in 22 days. This results in the yearly generation shares provided in Table \ref{tab:year}, which shows that the average share of renewable generation achieved by TA-UC (39.79\%) is slightly higher than the one yielded by CH-UC (39.64\%). It is worth mentioning that the cost reduction yielded by TA-UC does not imply an increase of the computational burden. In fact, the average time to solve the day-ahead CH-UC and TA-UC are 5 and 7 seconds, respectively. Both methods require computational times of the same order of magnitude since both consider the same number of time periods, variables, and constraints.

\renewcommand{\arraystretch}{1.2}
\setlength{\tabcolsep}{1.5mm}
\begin{table}
    \centering
    \caption{\textsc{Yearly Generation Shares for 2017}}
    \begin{tabular}{cccccc}
    \hline
Model & Wind (\%) & Solar (\%) & Base (\%) & Medium (\%) & Peak (\%) \\
    \hline
    CH-UC & 19.77 & 19.87 & 41.28 & 18.47 & 0.60 \\
    TA-UC & 19.85 & 19.94 & 41.27 & 18.37 & 0.56 \\
    \hline
    \end{tabular} 
    \label{tab:year}
\end{table}

As shown in the first part of this analysis, the benefits of the proposed model depend on both the penetration level and the type of renewable generation. To investigate this aspect further, Table \ref{tab:impact_renew} includes simulation results for the whole year 2017 under different penetration levels of wind and solar power production. Results in rows 2--6 are computed considering the same penetration level for wind and solar, while those in rows 7--11 and 12--16 refer to cases with wind and solar production only, respectively. This table shows the yearly cost savings achieved by the TA-UC model compared to the CH-UC, and the number of days for which the TA-UC provides lower, the same, or higher operating costs than the CH-UC, denoted as \# TA$<$CH, \# TA$=$CH, and \# TA$>$CH.

\renewcommand{\arraystretch}{1.2}
\setlength{\tabcolsep}{1.5mm}
\begin{table}
    \centering
    \caption{\textsc{Impact of Renewable Shares}}
    \begin{tabular}{cccccc}
    \hline
Wind (\%) & Solar (\%) & $\Delta C$ (\%) & \# TA$<$CH  & \# TA$=$CH & \# TA$>$CH \\
    \hline
    10 & 10 & 0.01 & 170 & 151 & 42 \\
    15 & 15 & 0.11 & 265 & 53 & 45 \\
    20 & 20 & 0.50 & 332 & 9 & 22 \\
    25 & 25 & 1.23 & 352 & 2 & 9 \\
    30 & 30 & 2.35 & 359 & 1 & 3 \\
    \hline
    20 & 0 & 0.01 & 233 & 85 & 45 \\
    30 & 0 & 0.08 & 290 & 32 & 41 \\
    40 & 0 & 0.21 & 322 & 12 & 29 \\
    50 & 0 & 0.37 & 320 & 10 & 33 \\
    60 & 0 & 0.56 & 317 & 24 & 22 \\
    \hline
    0 & 20 & 0.12 & 231 & 100 & 32 \\
    0 & 30 & 0.70 & 303 & 51 & 9 \\
    0 & 40 & 1.57 & 330 & 30 & 3 \\
    0 & 50 & 2.19 & 340 & 20 & 3 \\
    0 & 60 & 2.56 & 349 & 12 & 2 \\
    \hline
    \end{tabular} 
    \label{tab:impact_renew}
\end{table}

The main conclusion that can be derived from Table \ref{tab:impact_renew} is that the proposed model achieves significant yearly cost savings compared to the conventional approach in all cases without increasing the computational burden. Such cost savings range from 0.01\% to 2.56\% for the analyzed cases. By looking at any group of five rows of Table \ref{tab:impact_renew}, we can observe that an increase in the renewable penetration level leads to higher yearly costs savings. This is to be expected since renewable generation creates steep net demand variations that are poorly managed by the conventional CH-UC. Finally, if we compare cases with the same renewable penetration level but differently split between wind and solar, it can be concluded that the higher the renewable production coming from solar power plants, the larger the cost savings. This is so because solar generation has abrupt changes around sunrise and sunset, thus giving rise to a net demand curve with very steep variations, commonly known as the ``duck curve.'' As evidenced throughout this case study, the proposed TA-UC model outperforms the conventional CH-UC approach in the presence of sudden net demand variations. 

Another important factor that may affect the performance of the proposed method is the flexibility of the thermal generation portfolio. To analyze this aspect, we provide in Table \ref{tab:impact_thermal} the range of cost savings for nine different cases. Results obtained if the renewable generation comes from both wind and solar, only from wind or only from solar correspond to rows 2, 3 and 4, in that order. Results in column 2 are just a summary of those presented in Table \ref{tab:impact_renew}. The results for the \emph{high-flex} case (column 3) are obtained considering that $g1$--$g7$ are medium-load units, that is, their dispatch can be modified in real-time operation. Similarly, the results for the \emph{low-flex} case (column 4) are attained if $g1$--$g7$ are assumed as base-load units, that is, their dispatch must be fixed 24 hours in advance. As expected, if the thermal general portfolio is more inflexible, day-ahead dispatch and commitment decisions become more crucial to ensure an efficient real-time operation of the power system. For this reason, the proposed TA-UC achieves higher cost savings for the \emph{low-flex} case.

\renewcommand{\arraystretch}{1.2}
\setlength{\tabcolsep}{1.5mm}
\begin{table}
    \centering
    \caption{\textsc{Impact of Thermal Generation Flexibility}}
    \begin{tabular}{cccc}
    \hline
 & \emph{Base} case & \emph{High-flex} case & \emph{Low-flex} case\\
    \hline
Wind \& Solar & [0.01 -- 2.35] & [0.00 -- 1.04] & [0.27 -- 3.49] \\  
Only Wind & [0.01 -- 0.56] & [0.01 -- 0.08] & [0.30 -- 1.02] \\
Only Solar & [0.12 -- 2.56] & [0.07 -- 1.43] & [0.53 -- 4.76] \\
    \hline
    \end{tabular} 
    \label{tab:impact_thermal}
\end{table}

To conclude this analysis, the proposed method has been tested on a higher dimension case study inspired by the Spanish power system. The thermal generation portfolio of Table \ref{tab:gen_spain} includes 8 base-load units, 12 medium-load units, and 50 peak-load units. Some parameters have been slightly modified within a range to account for differences among units of the same type. The demand level and the capacity factor of wind and solar power generation correspond to year 2017 in Spain. The load shedding cost is  \euro10\,000/MWh.

\renewcommand{\arraystretch}{1.2}
\setlength{\tabcolsep}{0.9mm}
\begin{table}
\caption{\textsc{Thermal Generating Unit Data -- Spanish Power System}}
    \centering
    \begin{tabular}{cccccccc}
    \hline
 Type & $\underline{P}_g^G$ & $\overline{P}_g^G$ & $RU_g=RD_g$ & $UT_g=DT_g$ & $C^{SU}_{gt}$ & $C^M_{gt}$ & \# units \\
    \hline
    Base & 500 & 1000 & 3000 & [7.5--9] & 2000 & [20--25] & 8 \\
    Med & 400 & 800 & 4800 & [4--5] & 1000 & [50--55] & 12 \\
Peak & 0 & 500 & 4500 & - & 500 & [80--90] & 50 \\
    \hline
    \end{tabular} 
    \label{tab:gen_spain}
\end{table}

Results in Table \ref{tab:results_spain} show that TA-UC outperforms CH-UC in the three analyzed cases, all with a renewable share of 50\%. The highest cost saving amounts to 1.28\% and is achieved if all renewable generation comes from solar power units. Notice also that for most days of the year, the proposed method yields cheaper operating decisions than the conventional one. The computational time to solve the day-ahead UC problem for both approaches is lower than 10 seconds.

\renewcommand{\arraystretch}{1.2}
\setlength{\tabcolsep}{1.5mm}
\begin{table}
    \centering
    \caption{\textsc{Results -- Spanish Power System}}
    \begin{tabular}{cccccc}
    \hline
Wind (\%) & Solar (\%) & $\Delta C$ (\%) & \# TA$<$CH  & \# TA$=$CH & \# TA$>$CH \\
    \hline
    25 & 25 & 0.72 & 342 & 15 & 6 \\
    50 & 0 & 0.19 & 302 & 31 & 30 \\
    0 & 50 & 1.28 & 332 & 31 & 0 \\
    \hline
    \end{tabular} 
    \label{tab:results_spain}
\end{table}


\section{Conclusions} \label{sec:conclusions}


This paper proposes a novel time-adaptive unit commitment formulation to determine more efficiently the day-ahead status of thermal generating units. As done in the conventional unit commitment problem, our model also divides the decision horizon of one day into 24 time periods. However, instead of assuming time periods of one hour, we consider 24 time periods of different duration in order to better capture the varying net demand profile of the following day. Such durations are computed based on a chronological-time clustering technique.

As demonstrated by the numerical results, the proposed model leads to cost savings since it captures the intra-day variations of the net demand more precisely than the conventional approach. The yearly cost savings amount to 0.2--2.6\% for renewable penetration levels around 40--60\%. Yearly cost savings are particularly significant in power systems with a high solar penetration level. Such technology creates abrupt and monotonic net demand variations around sunrise and sunset (the duck curve), which are more accurately accounted for by the proposed approach. Finally, results also show that the lower the flexibility provided by the thermal generation, the higher the cost savings achieved by the proposed time-adaptive day-ahead unit commitment model.

Further research is required to investigate the performance of the proposed time-adaptive unit commitment via stochastic programming. Another ground for future research is the extension of the proposed approach to a rolling horizon framework under which the unit commitment is solved several times within a day as new information becomes available.



%





\ifCLASSOPTIONcaptionsoff
  \newpage
\fi

\balance



\bibliographystyle{IEEEtran}
\bibliography{mendeley_v2}
\end{document}